\documentclass[12pt]{article} 
\usepackage{amsmath, amssymb,hyperref,mathtools} 
\usepackage{amsthm}
\usepackage{enumerate}  
\usepackage{url} 

\usepackage{graphicx}

\newcommand\cx{{\mathbb C}}

\newcommand\ints{{\mathbb Z}}
\newcommand\re{{\mathbb R}}

\newcommand\cA{\mathcal{A}}
\newcommand\cB{\mathcal{B}}
\newcommand\cH{\mathcal{H}}

\DeclarePairedDelimiter\abs{\lvert}{\rvert}%
\newcommand\comp[1]{{\mkern2mu\overline{\mkern-2mu#1}}}
\newcommand\pmat[1]{\begin{pmatrix} #1 \end{pmatrix}}
\newcommand\seq[4]{#1_{#2},#1_{#3},\ldots,#1_{#4}}

\newtheoremstyle{plainsl}%
	{\topsep}
	{\topsep}
	{\slshape} 
	{}
	{\normalfont\bfseries}
	{.}
	{ }
	{}

\swapnumbers

{\theoremstyle{plainsl}
\newtheorem{theorem}{Theorem}[section]
\newtheorem{lemma}[theorem]{Lemma}
\newtheorem{corollary}[theorem]{Corollary}}
{\theoremstyle{remark}
}

\renewcommand\proof{\noindent\textsl{Proof. }}
\newcommand\sqr[2]{{\vbox{\hrule height.#2pt
    \hbox{\vrule width.#2pt height#1pt \kern#1pt
        \vrule width.#2pt}\hrule height.#2pt}}}
\renewcommand\qed{%
	\ifmmode\eqno\sqr53
	\else\nolinebreak\ \hfill\sqr53\medbreak\fi}

\DeclareMathOperator{\wt}{wt}

\newcommand\one{{\bf1}}
\newcommand\zero{{\bf0}}

\newcommand\grp[1]{\langle #1\rangle}

\title{Uniform Mixing on Cayley Graphs}
\author{Chris Godsil, Hanmeng Zhan}

\begin{document}
\maketitle

\begin{abstract}
We provide new examples of Cayley graphs on which the quantum walks reach uniform mixing. Our first result is a complete characterization of all $2(d+2)$-regular Cayley graphs over $\mathbb{Z}_3^d$ that admit uniform mixing at time $2\pi/9$. Our second result shows that for every integer $k\ge 3$, we can construct Cayley graphs over $\mathbb{Z}_q^d$ that admit uniform mixing at time $2\pi/q^k$, where $q=3, 4$.

We also find the first family of irregular graphs, the Cartesian powers of the star $K_{1,3}$, that admit uniform mixing. 

\vspace{.1in}
\par\noindent{\em Keywords}: 
quantum walk, uniform mixing, Cayley graph
\end{abstract}

\section{Introduction \label{intro}}

A continuous-time quantum walk on a graph $X$ is defined by the transition matrix
\[U(t):=\exp(itA)=\sum_{k\ge 0}\frac{(itA)^k}{k!},\]
where $A$ is the adjacency matrix of $X$.  The probability that at time $t$, the quantum walk with initial state represented by $u$ is in the state represented by $v$ is
\[\left| U(t)_{uv}\right|^2.\]
We say $X$ admits \textsl{uniform mixing} at time $t$ if the above probability is the same for all vertices $u$ and $v$. A weaker version of uniform mixing, called \textsl{local uniform} mixing, occurs if the probability distribution given by a column of $U(t)$ is uniform.

Uniform mixing on graphs is rare, and almost all the known examples are Cayley graphs. The current list of Cayley graphs contains the complete graphs $K_2$, $K_3$ and $K_4$ \cite{Ahmadi2003}, the Hamming graphs $H(d,2)$, $H(d,3)$ and $H(d,4)$ \cite{Carlson2006,Moore2001}, the Paley graph of order nine \cite{Godsil2013}, some strongly regular graphs from regular symmetric Hadamard matrices \cite{Godsil2013}, some linear Cayley graphs over $\ints_2^d$, $\ints_3^d$ and $\ints_4^d$ \cite{Chan2013,Mullin2013}, and the Cartesian product of graphs which admit uniform mixing at the same time. These graphs share the following features.

\begin{enumerate}[(a)]
\item
They are Cayley graphs over abelian groups;
\item
They have integer eigenvalues;
\item
Their mixing times are rational multiples of $\pi$.
\end{enumerate}

That being said, a graph that admits uniform mixing is not necessarily vertex-transitive or even regular. As we will see in Section \ref{star}, the star $K_{1,3}$ admits uniform mixing at time $2\pi/\sqrt{27}$. The Cartesian powers of $K_{1,3}$ are so far the only family of irregular graphs found to admit uniform mixing, and the only family to which none of the above features applies.

Although characterizing uniform mixing in general seems daunting, the problem can be reduced if we require some regularity on the graphs. The vertex transitive graphs, which provide most of the known examples, satisfy the property that for any two vertices $u$ and $v$, there is a graph automorphism that maps $u$ to $v$. The transition matrix of a vertex-transitive graph is entirely determined by one of its rows. Therefore for all Cayley graphs, uniform mixing is equivalent to local uniform mixing. 

In this paper, we provide new examples of Cayley graphs over $\ints_q^d$ that admit uniform mixing. Our examples contain both an infinite family of Cayley graphs on which mixing occurs at time $2\pi/9$, and infinite families of Cayley graphs on which mixing occurs arbitrarily faster. Theorem \ref{2-gen} provides a complete characterization of $2(d+2)$-regular Cayley graphs over $\ints_3^d$ that admit uniform mixing at time $2\pi/9$. Theorem \ref{q=3(1)}, Theorem \ref{q=4} and Theorem \ref{q=3(2)} show that, for an arbitrarily large integer $k$, we can construct families of Cayley graphs over $\ints_q^d$ that admit uniform mixing at time $2\pi/q^k$, where $q=3, 4$. These examples extend the results of Mullin \cite{Mullin2013} and Chan \cite{Chan2013}.

\section{Quotients of Hamming Graphs}

A \textsl{Cayley graph} over the additive group $\ints_q^d$ is a graph with vertex set $\ints_q^d$ and edge set
\[\{(g,h):h-g\in C\}\]
for some subset $C$ of $\ints_q^d$. We will follow Godsil and Royle \cite{Godsil2001} and denote this graph by $X(\ints_q^d, C)$, and call $C$ the \textsl{connection set}. Further, we assume $C$ is inverse-closed and does not contain the identity element, so that $X(\ints_q^d,C)$ is a simple graph. Note here $q$ is not necessarily a prime power -- we will think of the vertex set $\ints_q^d$ as a module in this paper.

While no extra algebraic structure of $\ints_q^d$ is needed to define a Cayley graph, it is often convenient to view the underlying group $\ints_q^d$ as a $\ints_q$-module. Two easy observations follow. First, if $\phi$ is a module automorphism of $\ints_q^d$, then the two Cayley graphs $X(\ints_q^d, C)$ and $X(\ints_q^d,\phi(C))$ are isomorphic.
Second, $X(\ints_q^d, C)$ is connected if and only if its connection set contains a basis of $\ints_q^d$. An example of connected Cayley graphs over $\ints_q^d$ is the \textsl{Hamming graph} $H(d,q)$, whose connection set consists of non-zero multiples of the standard basis $\{\seq{e}{1}{2}{d}\}$ of $\ints_q^d$.

We will pay special attention to the quotient graphs of Hamming graphs, as they form an important family of Cayley graphs over $\ints_q^d$. The \textsl{Hamming distance} of two elements in $\ints_q^d$ is the number of coordinates in which they differ. Consider a submodule $\Gamma$ of $\ints_q^d$ with Hamming distance at least three. The partition of $\ints_q^d$ by the cosets of $\Gamma$ satisfies three properties: 
\begin{enumerate}[(i)]
\item every coset is a coclique;
\item every vertex is adjacent to at most one vertex in a coset;
\item if some vertex in the coset $g+\Gamma$ is adjacent to a vertex in another coset $h+\Gamma$, then there is a matching between $g+\Gamma$ and $h+\Gamma$.
\end{enumerate}
The \textsl{quotient graph} of $H(d,q)$ induced by $\Gamma$, denoted $H(d,q)/\Gamma$, is a graph with a vertex for each coset of $\Gamma$, such that two vertices are adjacent if there is a matching between the two associated cosets. We note that every quotient graph of $H(d,q)$ is a Cayley graph for a quotient module of $\ints_q^d$.

\begin{lemma}
Let $C$ be the connection set of the Hamming graph $H(d,q)$. For a submodule $\Gamma$ of $\ints_q^d$ with Hamming distance at least three, let 
\[C/\Gamma=\{c+\Gamma:c\in C\}.\]
Then
\[H(d,q)/\Gamma\cong  X(\ints_q^d/\Gamma, C/\Gamma).\]
\label{qmodule}
\end{lemma}

\proof
Let $g+\Gamma$ and $h+\Gamma$ be two vertices of $X(\ints_q^d,C)/\Gamma$. They are adjacent if and only if there exist $x,y\in\Gamma$ and $c\in C$ such that $(g+x)-(h+y)=c$, or equivalently, \[(g+\Gamma)-(h+\Gamma)=c+\Gamma\in C+\Gamma. \tag*{\sqr53}\]

\section{Linear Cayley Graphs}
A Hamming graph $H(d,q)$ is known to admit uniform mixing if and only if $q\in\{2,3,4\}$. As uniform mixing on $H(d,q)$ implies uniform mixing on some of its quotient graphs, for our interest it is important to understand which Cayley graphs are quotients of $H(d,q)$. In this section, we show that quotient graphs of $H(d,q)$ are exactly  the \textsl{linear} Cayley graphs over $\ints_q^d$, that is, graphs $X(\ints_q^d, C)$ for which $C\cup\{\mathbf{0}\}$ is closed under multiplication by $\ints_q$.

Since $\ints_2$ has only one non-zero element, the connection set of every Cayley graph over $\ints_2^d$ is trivially closed under multiplication of the non-zero elements of $\ints_2$. Similarly, for Cayley graphs over $\ints_3^d$, since we assume the connection set is inverse-closed, it is closed under multiplication of the only two non-zero elements of $\ints_3$. Thus, all Cayley graphs over $\ints_2^d$ or $\ints_3^d$ are linear, and we can characterize uniform mixing on them in terms of the submodules that induce the quotients.

\begin{theorem}
Let $X(\ints_q^r, C)$ be a connected linear Cayley graph with valency $d(q-1)$. Then $X(\ints_q^r, C)$ is isomorphic to a quotient graph $H(d,q)/\Gamma$ for some submodule $\Gamma$ of $\ints_q^d$, where $|\Gamma|=q^{d-r}$ and $\Gamma$ has Hamming distance at least three. 
\end{theorem}

\proof
Let $C$ be the connection set of $X$. We can partition $C$ into cells $\seq{C}{1}{2}{d}$ such that two elements lie in the same cell if and only if they are multiples of each other, and the first $r$ cells contains a basis of $\ints_q^d$. Since $C_j$ is a cyclic group, we may assume $C_j=\left<v_j\right>$. Define a module homomorphism from $C$ to $\ints_q^r\times \ints_q^{d-r}$ by
\[f(v_j)=\begin{cases}(v_j,0), \quad\text{ if }v\in C_j \text{ for } j\in\{1,2,\cdots, r\},\\(v_j,0) + e_j,\quad \text{ if } v_j\in C_j \text{ for } j\in\{r+1,r+2,\cdots,d\},\end{cases}\]
where $\{\seq{e}{1}{2}{d}\}$ is the standard basis of $\ints_q^d$. Let
\[ C'=\{f(v):v\in C\}.\]
Then $ C'$ consists of all non-zero multiples of a basis of $\ints_q^r\times\ints_q^{d-r}$. Thus, $X(\ints_q^d, C')$ is isomorphic to the Hamming graph $H(d,q)$. Let $\phi$ be a module automorphism of $\ints_q^d$ that maps a basis in $C'$ to the standard basis. Let $\Gamma'$ be the submodule of $\ints_q^d$ generated by $\{\seq{e}{r+1}{r+2}{d}\}$, and  $\Gamma$ the image of $\Gamma'$ under $\phi$. Clearly $|\Gamma|=q^{d-r}$. By Lemma \ref{qmodule},
\[H(d,q)/\Gamma\cong X(\ints_q^d, C')/\Gamma' \cong X(\ints_q^r, C). \]
Since we started with a simple Cayley graph, $\Gamma$ must have Hamming distance at least three.
\qed

Conversely, for any quotient graph $H(d,q)/\Gamma$, where $\Gamma$ has minimum distance at least three, we can find a connection set of the linear Cayley graph isomorphic to $H(d,q)/\Gamma$.

\begin{theorem}
Let $\Gamma$ be a submodule of $\ints_q^d$ with size $q^{d-r}$ and Hamming distance at least three. Let $Q$ be a parity check matrix of $\Gamma$, and $C$ the set of non-zero multiples of the columns of $Q$. Then $H(d,q)/\Gamma$ is isomorphic to the $d(q-1)$-regular graph $X(\ints_q^r, C)$.
\label{cnset}
\end{theorem}

\proof
Without loss of generality, we may assume $\Gamma$ is generated by the columns of the following block matrix
\[\pmat{R\\S},\]
where $S$ is square and invertible over $\ints_q$. Let
\[P=\pmat{I & -RS^{-1}\\0 & S^{-1}}.\]
We have
\[P\pmat{R\\S}=\pmat{0\\I},\]
that is, $P$ defines a module automorphism $\phi$ of $\ints_q^d$ that maps $\Gamma$ to the submodule generated by $\{\seq{e}{r+1}{r+2}{d}\}$. 
Note that the partitioned matrix
\[Q=\pmat{I & -RS^{-1}}\]
is a parity check matrix of $\Gamma$. If $D$ is the connection set of $H(d,q)$ and $C$ the multiples of columns of $Q$, then by Lemma \ref{qmodule},
\begin{align*}
H(d,q)/\Gamma &\cong X(\ints_q^d/\Gamma, D/\Gamma)\\
& \cong X(\ints_q^d/\phi(\Gamma), \phi(D)/\phi(\Gamma))\\
& \cong X(\ints_q^d, C).
\end{align*}
Finally, since the minimum distance of $\Gamma$ is the minimum number of linearly dependent columns of its parity check matrix $Q$, no two columns of $Q$ are multiples of each other. It follows that $X(\ints_q^r, C)$ has valency $d(q-1)$.
\qed

The quotient approach provides a convenient way to characterize uniform mixing on linear Cayley graphs. If $X = H(d,q)/\Gamma$, the vertices of $X$ are cosets $\Gamma + v$ of the subgroup, and each entry of $U_X(t)$ is a block sum of the entries of $U_{H(d,q)}(t)$. More specifically, the $(0, v)$-entry of $U_X(t)$ can be expressed as follows. For more details, see Mullin's thesis \cite[Ch 8]{Mullin2013}.

\begin{theorem}[Mullin]
Let $X = H(d,q)/\Gamma$. We have
\[U(t)_{0, v}=\left(\frac{e^{-it}}{q}\right)^d \sum_{a \in \Gamma +v} (e^{qit}+q-1)^{d-\wt(a)} (e^{qit}-1)^{\wt(a)}. \tag*{\sqr53}\]
\label{Mullin}
\end{theorem}

For each coset $\Gamma+v$, let $W_v(x,y)$ denote its homogeneous weight enumerator. Note that the right hand side of the above equation is 
\[W_v(e^{qit}+q-1, e^{qit}-1).\]
For the $(0,0)$-entry, MacWilliams' identity simplifies the expression to
\[U(t)_{0,0} = \left(\frac{e^{-it}}{q}\right)^d \abs{\Gamma} W_{\Gamma^{\perp}}(e^{qit},1).\]
Thus we have a necessary condition for uniform mixing.

\begin{corollary}
If $H(d,q)/\Gamma$ admits uniform mixing at time $t$, then
\[\abs{W_{\Gamma^{\perp}}(e^{qit}, 1)}^2 = \abs{\Gamma^{\perp}}. \tag*{\sqr53}\]
\end{corollary}

We already know that uniform mixing occurs on $H(d,2), H(d,3)$ and $H(d,4)$ at time $\pi/4, 2\pi/9$ and $\pi/4$, respectively. It is natural to see if their quotients also admit uniform mixing at these special times. 

\begin{corollary}
\begin{enumerate}[(a)]
\item $H(d,2)/\Gamma$ admits uniform mixing at time $\pi/4$ if and only if for each coset $\Gamma+v$,
\[\abs{W_v(i,1)}^2 = \abs{\Gamma};\]
\item $H(d,3)/\Gamma$ admits uniform mixing at time $2\pi/9$ if and only if for each coset $\Gamma+v$,
\[\abs{W_v(e^{2\pi i/3},1)}^2 = \abs{\Gamma};\]
\item $H(d,4)/\Gamma$ admits uniform mixing at time $\pi/4$ if and only if for each coset $\Gamma+v$,
\[\abs{W_v (-1,1)}^2 = \abs{\Gamma}.\]
\end{enumerate}
\end{corollary}

Consider the binary $[17,9,5]$-quadratic code. A numerical check on the weight distributions of its cosets shows that it admits uniform mixing at time $\pi/4$. 
%
%
%
%

\section{Quotient Graphs with One Generator \label{quotient}}

For a quotient graph $H(d,q)/\Gamma$, the entries of its transition matrix are block sums of the transition matrix of $H(d,q)$. As functions of the time $t$, these block sums can be greatly simplified if we plug in the time $\tau_q$ when $H(d,q)$ admits uniform mixing. Thus, at the specific time $\tau_q$, whether uniform mixing occurs on $H(d,q)/\Gamma$ is fully determined by the weight distributions of the cosets of $\Gamma$. For details see Mullin's Ph.D. thesis \cite[Ch 8]{Mullin2013}.

The \textsl{Hamming weight} $\wt(a)$ of an element $a$ in $\ints_q^d$ is the number of non-zero entries of $a$. In this section, we characterize quotients $H(d,q)/\grp{a}$ that admit uniform mixing at time $\tau_q$, where the submodule is generated by one element $a$ with Hamming weight $\wt(a)$ at least three. As a special case, the ``folded" Hamming graphs $H(d,q)/\grp{\one}$ have been studied in \cite{Mullin2013}. In general, Theorem \ref{cnset} gives the matrix $Q$ representing the connection set of $H(d,q)/\grp{a}$. By row reduction, column permutation, and column scaling of $Q$, it is not hard to see that $H(d,q)/\grp{a}$ is a Cartesian product of a Hamming graph and a folded Hamming graph. Combining this observation and the results on the folded Hamming graphs, we give the following characterization.

\begin{theorem}
Let $a$ be a vector in $\ints_q^d$ with $\wt(a)\ge 3$, where $q=2,3,4$. Then
\begin{enumerate}[(a)]
\item
Uniform mixing occurs on $H(d,2)/\grp{a}$ at time $\tau_2=\pi/4$ if and only if $\wt(a)$ is odd;
\item
Uniform mixing occurs on $H(d,3)/\grp{a}$ at time $\tau_3=2\pi/9$ if and only if $\wt(a)$ is not divisible by three;
\item 
Uniform mixing occurs on $H(d,4)/\grp{a}$ at time $\tau_4=\pi/4$ if and only if $\wt(a)$ is odd.
\qed
\end{enumerate}
\end{theorem}
This theorem takes care of all the connected $(d+1)(q-1)$-regular Cayley graphs over $\ints_q^d$, which admit uniform mixing at $\tau_q$, for $q\in\{2,3,4\}$.

\section{Quotient Graphs with Two Generators} 
We move on to the quotients of Hamming graphs where the submodules are generated by two elements. Our goal is to characterize the $2(d+2)$-regular Cayley graphs over $\ints_3^d$ that admit uniform mixing at time $2\pi/9$, as given in Theorem \ref{2-gen}. The tool we use is the \textsl{weight distribution} of a code, that is, the sequence
\[(\alpha_0, \alpha_1, \alpha_2, \cdots )\]
where $\alpha_j$ is the number of of codewords with Hamming weight $j$. To begin, we need the following conditions on the cosets of $\Gamma$ from \cite[Ch 8]{Mullin2013}.

\begin{theorem}[Mullin]
Let $\Gamma$ be a submodule of $\ints_3^d$ with Hamming distance at least three such that $|\Gamma|=3^s$. For any coset of $\Gamma$, let $n_j$ be the number of elements in it with weight $j$ modulo three. Uniform mixing occurs on $H(d,q)/\Gamma$ at time $2\pi/9$ if and only if the weight distribution of every coset of $\Gamma$ satisfies 
\[n_0n_1+n_0n_2+n_1n_2=3^{2s-1}-3^{s-1}.\tag*{\sqr53}\]
\label{q3}
\end{theorem}

Since we have to examine the weight distribution of every coset of $\Gamma$, it helps to understand the relation between the weights of $\Gamma$ and the weights of $\Gamma+c$ for each vector $c$. Suppose $\Gamma$ is generated by $s$ elements, and let $M$ be the matrix with the generators of $\Gamma$ as its columns. Then each element in $\Gamma$ can be written as $My$ for some vector $y\in \ints_3^s$, and is uniquely associated to an element $My+c$ in the coset $\Gamma+c$. We will refer to 
\[\wt\left(My+c\right)-\wt\left(My\right)\]
as the \textsl{weight change} of the element $My$ with respect to $c$. The following gives the necessary condition on the weight changes for both $\Gamma$ and $\Gamma+c$ to satisfy the weight distribution condition in Theorem \ref{q3}.

\begin{lemma}
Let $\Gamma$ be a submodule of $\ints_3^d$ with size $3^s$ and minimum distance three, and let $\Gamma+c$ be a coset of $\Gamma$. For $j=0,1,2$, let $\Gamma_j$ denote the set of elements in $\Gamma$ with weight congruent to $j$ modulo three. Let $m_j$ be the number of elements in $\Gamma_j$ whose weight change with respect to $c$ is one modulo three. If both $\Gamma$ and $\Gamma+c$ satisfy the weight distribution condition in Theorem \ref{q3}, then either $m_j=0$ for all $j$, or $m_j$ satisfies the following:
\begin{gather*}
m_0+m_1+m_2=3^{s-1}\\
(m_0n_0+m_1n_1+m_2n_2)+3(m_0m_1+m_1m_2+m_0m_2)=3^{2s-1}-3^{2s-2}.
\end{gather*}
\label{wtchange}
\end{lemma}
\proof
We will calculate the weights over $\ints_3$. Let $My$ be an element in $\Gamma$. Notice that
\[\wt(My+c)=(My+c)^T(My+c)=\wt(My)+\wt(c)+2c^TMy.\]
Since $\wt(c)$ depends only on $c$, and the condition in Theorem \ref{q3} is symmetric on $n_0,n_1,n_2$, we may assume without loss of generality that $\wt(c)=0$. If $c^TM=0$, then the weight change of each element in $\Gamma$ is zero. Otherwise, there are exactly $3^{s-1}$ vectors $y$ such that 
\[2c^TMy=1.\]
Therefore 
\[m_0+m_1+m_2=3^{s-1}.\]
Notice that for each solution $y$ to $2c^TMy=1$, the vector $2y$ is a solution to $2c^TMy=2$. Thus there are equal number of elements in $\Gamma_j$ with weight change one and weight change two. It follows that in the coset $\Gamma+c$, the number of elements with weight $j$ is 
\[n'_j=n_j-2m_j+m_{j-1}+m_{j+1},\]
where the subcripts are calculated modulo three. Since $\Gamma+c$ satisfies the weight distribution conditon in Theorem \ref{q3}, 
\[n'_0n'_1+n'_0n'_2+n'_1n'_2=3^{2s-1}-3^{s-1}.\]
This together with the fact that
\[n_0n_1+n_0n_2+n_1n_2=3^{2s-1}-3^{s-1}\]
yields
\[(m_0n_0+m_1n_1+m_2n_2)+3(m_0m_1+m_1m_2+m_0m_2)=3^{2s-1}-3^{2s-2}.\tag*{\sqr53}\]

With the above observation, we characterize the quotient graphs $H(d,q)/\grp{a,b}$ that admit uniform mixing at time $2\pi/9$ in terms of the generators $a$ and $b$. 

\begin{theorem}
Uniform mixing occurs on $H(d,3)/\grp{a,b}$ at time $2\pi/9$ if and only if one of the following holds:
	\begin{enumerate}[(i)]
	\item $a^Tb\equiv 0\pmod{3}$, $\wt(a)\not\equiv 0\pmod{3}$, and $\wt(b)\not \equiv 0\pmod{3}$,
	\item $a^Tb\not\equiv 0\pmod{3}$, and $\wt(a)\not\equiv \wt(b)\pmod{3}$ unless $\wt(a)\equiv \wt(b)\equiv 0\pmod{3}$.
	\end{enumerate}
\label{2-gen}
\end{theorem}

\proof
For notational convenience, we define the \textsl{weight structure} of the coset $\Gamma+c$ to be the tuple with coordinates $n_0,n_1,n_2$ in non-descending order, denoted by $W(\Gamma+c)$. By Theorem \ref{q3}, the quotient graph $H(d,3)/\grp{a,b}$ admits uniform mixing at time $2\pi/9$ if and only if the weight distribution of every coset $\Gamma+c$ satisfies 
\begin{gather*}
n_0n_1+n_0n_2+n_1n_2=24\\
n_0+n_1+n_2=9
\end{gather*}
which holds if and only if for all $c$,
\[W(\Gamma+c)\in\{(1,4,4),(2,2,5)\}.\]

We first show that $W(\Gamma)$ lies in the above set if and only if one of the conditions (i) and (ii) holds. Let 
\[M=\pmat{a&b}\]
be a matrix and let 
\[y=\pmat{y_1\\y_2}\]
be a vector in $\ints_3^2$. The weight of $My$ is 
\[\wt(My)=y^TM^TMy=\wt(y_1)\wt(a)+\wt(y_2)\wt(b)+2y_1y_2a^Tb.\]
Thus the weights of the elements in $\Gamma$ are
\begin{center}
\begin{tabular}{|c|c|c|}
\hline
label & weight & multiplicity\\
\hline
$w_0$ & $0$ & 1\\
\hline
$w_1$ &$\wt(a)$ & 2\\
\hline 
$w_2$ & $\wt(b)$ & 2\\
\hline
$w_3$ & $\wt(a)+\wt(b)+a^Tb$ & 2\\
\hline
$w_4$ & $\wt(a)+\wt(b)+2a^Tb$ & 2\\
\hline
\end{tabular}
\end{center}
Since $\Gamma$ is a group of order nine, $n_0$ is odd and $n_1,n_2$ are even. We consider two cases.
\begin{enumerate}[(a)]
\item Suppose 
\[W(\Gamma)=(1,4,4).\]
Then half of $\{w_1,w_2,w_3,w_4\}$ are one, and the rest are two.
\begin{itemize}
\item $\wt(a)=\wt(b)\ne 0$. Then $w_3=w_4\notin\{0, w_1\}$ if and only if $a^Tb=0$.
\item $\wt(a)=2\wt(b)\ne 0$. Then $w_3=a^Tb=2w_4$. It follows that two of $\{w_1,w_2,w_3,w_4\}$ are one and the others are two if and only if $a^Tb\ne 0$.
\end{itemize}
\item Suppose
\[W(\Gamma)=(2,2,5).\]
Then half of $\{w_1,w_2,w_3,w_4\}$ are zero, and the rest are one and two respectively. 
\begin{itemize}
\item $\wt(a)=\wt(b)=0$. Then $w_3=a^Tb=2w_4$. Thus $\{w_3,w_4\}=\{1,2\}$ if and only if $a^Tb\ne 0$.
\item $\wt(a)=2\wt(b)\ne 0$. Then $w_3=a^Tb=2w_4$. It follows that $w_3=w_4=0$ if and only if $a^Tb=0$.
\end{itemize}
\end{enumerate}
Summarizing the above yields the conditions (i) and (ii).

Next we show that if 
\[W(\Gamma)\in\{(1,4,4),(2,2,5)\}\]
then 
\[W(\Gamma+c)=W(\Gamma),\]
for all $c\in\ints_3^d$. By Lemma \ref{wtchange}, the weight changes $m_0,m_1,m_2$ are either zero, or satisfy
\begin{gather}
m_0+m_1+m_2=3, \notag \\
(m_0n_0+m_1n_1+m_2n_2)+3(m_0m_1+m_0m_2+m_0m)_3=18. \label{wtch}
\end{gather}
In the latter case, the weight changes are either
\[m_0=m_1=m_2=1\]
or
\[\{m_0,m_1,m_2\}=\{0,1,2\}.\]
It is easy to see that if $m_0=m_1=m_2=1$, then $W(\Gamma+c)=W(\Gamma)$. Suppose 
\[\{m_0,m_1,m_2\}=\{0,1,2\}.\]
Then by Equation (\ref{wtch}), when $W(\Gamma)=(1,4,4)$, we have 
\begin{gather*}
n_0=1,\quad n_1=n_2=4\\
m_0=0,\quad \{m_1,m_2\}=\{1,2\}
\end{gather*}
and when $W(\Gamma)=(2,2,5)$, we have
\begin{gather*}
n_0=5,\quad n_1=n_2=2\\
m_0=2,\quad \{m_1,m_2\}=\{0,1\}.
\end{gather*}

Since the weight distribution of $\Gamma+c$ is
\[n'_j=n_j-2m_j+m_{j-1}+m_{j+1}\]
again we have $W(\Gamma+c)=W(\Gamma)$.
\qed

It is perhaps surprising that in all the characterizations discussed so far, the condition for a quotient graph $H(d,q)/\Gamma$ to admit uniform mixing only relies on the group generators, although Theorem \ref{q3} suggests checking the weight distribution of every coset of this group. It would reduce the problem of checking uniform mixing on $H(d,q)/\Gamma$ at time $\tau_q$ considerably if this was true in general.

\section{Hamming Schemes}
In this section, we construct Cayley graphs over $\ints_q^d$ that admit uniform mixing earlier than the complete graph $K_q$. The construction is based on the association scheme which spans the adjacency algebra of a Hamming graph, called the Hamming scheme. We introduce the basic concepts of association schemes and some useful results on the eigenvalues of Hamming schemes.

An \textsl{association scheme} with $d$ classes is a set $\cA=\{\seq{A}{0}{1}{d}\}$ of 01-matrices that satisfies the following conditions:
\begin{itemize}
\item $A_0=I$.
\item $\sum_{r=0}^dA_j=J$.
\item $A_r^T\in \cA$ for $r=0,1,\cdots,n$.
\item $A_rA_s=A_rA_s\in\mathrm{span}(\cA)$.
\end{itemize}

The association scheme $\cA$ generates an algebra over $\cx$, which is referred to as the \textsl{Bose-Mesner algebra} of $\cA$. This algebra has an orthogonal basis of idempotents $\seq{E}{0}{1}{d}$. Thus for each matrix $A_r$ in the scheme, there are scalars $p_r^{(d)}(s)$ such that
\[A_r=\sum_{s=0}^d p_r^{(d)}(s)E_s.\]
These scalars are called the eigenvalues of the scheme $\cA$. When $d$ is clear from the context, we drop the superscript and write $p_r(s)$. The following theorem due to Chan \cite{Chan2013} shows that whether a graph admits uniform mixing depends only on its spectrum and the eigenvalues of the Bose-Mesner algebra containing its adjacency matrix.

\begin{theorem}[Chan]
Let $X$ be a graph on $v$ vertices whose adjacency matrix belongs to the Bose-Mesner algebra of $\cA$. Let $p_j(s)$ be the eigenvalues of $\cA$. Suppose the spectral decomposition of $A(X)$ is
\[A(X)=\sum_{s=0}^d\theta_sE_s.\]
The continuous quantum walk of $X$ is uniform mixing at time $\tau_q$ if and only if there exist scalars $\seq{t}{0}{1}{d}$ such that 
\begin{itemize}
\item $|t_0|=|t_1|=\cdots=|t_d|=1$,
\item $\sqrt{v}e^{i\tau_q\theta_s}=\sum_{j=0}^dp_j(s)t_j$ for $s=0,1,\cdots, d$.
\qed \label{coeff}
\end{itemize}
\end{theorem}

A \textsl{Hamming scheme} $\cH(d,q)$ is an association scheme constructed from the Hamming graph $H(d,q)$. The matrix $A_j$ is the adjacency matrix of the $r$-th distance graph of $H(d,q)$, which has the same vertex set as $H(d,q)$ such that two vertices are adjacent if they are at distance $r$ in $H(d,q)$. The eigenvalues of the scheme $\cH(d,q)$ satisfy
\begin{equation}
p^{(d)}_r(s)=[x^r](1+(q-1)x)^{d-s}(1-x)^s, \label{hamev}
\end{equation}
for $s=0,1,\cdots,d$. These are called the \textsl{Krawtchouk polynomials}. They satisfy the following properties.

\begin{lemma}
Let $p_r(s)$ be the eigenvalues of the Hamming scheme $\cH(d,q)$. Then
\begin{enumerate}[(i)]
\item $p_r(s)=\sum_h(-q)^h(q-1)^{r-h}\dbinom{d-h}{r-h}\dbinom{s}{h}$.
\item $p_r(s)-p_r(s-1)+(q-1)p_{r-1}(s)+p_{r-1}(s-1)=0$.
\item $p^{(d+1)}_r(s)-p^{(d+1)}_r(s+1)=qp^{(d)}_{r-1}(s)$.
\item If $q=2$, then 
\[p_{r-1}(s)-p_{r-1}(s+2)=4\sum_h (-2)^h\binom{d-2-h}{r-2-h}\binom{s}{h}.\]
\end{enumerate}
\label{properties}
\end{lemma}
\proof
By Equation (\ref{hamev}), 
\begin{align*}
p_r(s)&=[x^r](1+(q-1)x)^{d-s}(1+(q-1)x-qx)^s\\
&=[x^r]\sum_h \binom{s}{h}(1+(q-1)x)^{d-h}(-qx)^h\\
&=\sum_h \left([x^h]\binom{s}{h}(-qx)^h\right)\left([x^{r-h}](1+(q-1)x)^{d-h}\right)\\
&=\sum_h (-q)^h(q-1)^{r-h}\binom{d-h}{r-h}\binom{s}{h}. 
\end{align*}
Properties (ii) and (iii) follow from Equation (\ref{hamev}), and Property (iv) follows from the above three properties for $q=2$.
\qed

We will use property (iii) frequently, so we rephrase it in the matrix form. Let $P^{(d)}$ be the eigenvalue matrix for the scheme $\cH(d,q)$. Let 
\[C_d = \pmat{
	0&1&0&\cdots&0\\
	0&0&1&\cdots&0\\
	\cdots & \cdots & \cdots & \cdots & \cdots\\
	0&0&0&\cdots&1\\
	1&0&0&\cdots&0}\]
be an $d\times d$ circulant matrix. Then Property (iii) is equivalent to
\begin{equation}
\pmat{I_{d-1} & 0\\ 0 & 0} (I_d-C) P^{(d)} \pmat{I_{d-1} & 0\\ 0 & 0}= q P^{(d-1)}
\label{ev_rec}
\end{equation}

\section{Sufficient Conditions for Uniform Mixing in Hamming Schemes}
A graph on $n$ vertices admits uniform mixing at time $t$ if and only if $\sqrt{n}U(t)$ is equal to a complex Hadamard matrix. For notational convenience, let $\cB(q)$ denote the Bose-Mesner algebra of the Hamming scheme $\cH(d,q)$. For $q\in\{2,3,4\}$, we can construct complex Hadamard matrices in $\cB(q)$ from a primitive $q$-th roots of unity.


\begin{lemma}
Let $\zeta_q$ be a primitive $q$-th root of unity. For $q\in\{2,3,4\}$, the matrix
\[e^{i\beta}(I_q+\zeta_{6-q} (J_q-I_q))^{\otimes d}\]
is a complex Hadamard matrix in $\cB(q)$.
\qed
\end{lemma}
\proof
First note that $I_q+\zeta_{6-q}(J_q-I_q)$ is a complex Hadamard matrix of order $q$ for $q\in\{2,3,4\}$. In fact, it is a scalar multiple of $U_{K_q}(\tau_q)$ for some mixing time $\tau_q$ of $K_q$. Now if $H_1$ and $H_2$ are both complex Hadamard matrices of order $q$, then $H_1\otimes H_2$ is flat and thus a complex Hadamard matrix of order $q^2$. Lastly, a unimodular scalar multiple of a complex Hadamard matrix is again a complex Hadamard matrix.
\qed

The following generalizes Lemma 3.2 in \cite{Chan2013}.

\begin{lemma}
Let $X$ be a graph in $\cB(q)$ with eigenvalues $\seq{\theta}{0}{1}{d}$, and let $\epsilon\in\{1,-1\}$. Suppose $k\ge 2$.
\begin{enumerate}[(i)]
\item If $q=2$, and
\[\theta_s-\theta_0\equiv  \epsilon 2^{k-1}s\pmod{2^{k+1}},\]
for $s=0,1,\cdots,d$, then $X$ admits uniform mixing at time $\pi/2^k$.
\item If $q=3$, and 
\[\theta_s-\theta_0\equiv  \epsilon 3^{k-1}s \pmod{3^k},\]
for $s=0,1,\cdots,d$, then $X$ admits uniform mixing at time $2\pi/3^k$.
\item If $q=4$, and 
\[\theta_s-\theta_0\equiv 2^ks\pmod{2^{k+1}},\]
for $s=0,1,\cdots,d$, then $X$ admits uniform mixing at time $\pi/2^k$.
\end{enumerate}
\label{theta_diff}
\end{lemma}
\proof
Suppose $q\in\{2,3,4\}$. Let 
\begin{align*}
H_q&=e^{i\beta}(I_q+\zeta_{6-q} (J_q-I_q))^{\otimes d}\\
&=e^{i\beta}\left((1+(q-1)\zeta_{6-q})\left(\frac{1}{q}J_q\right)+(1-\zeta_{6-q})\left(I_q-\frac{1}{q}J_q\right)\right)^{\otimes d}.
\end{align*}
Suppose $A\in\cB(q)$ has spectral decomposition
\[A=\sum_{r=0}^d \theta_r E_r.\]
By Equation (4.2) in \cite[Sec 4]{Godsil2010a}, 
\[H_q=e^{i\beta}\sum_{r=0}^d(\zeta_{6-q})^rA_r,\]
where $A_r$ is the adjacency matrix of the $r$-distance graph of $H(d,q)$. Hence the condition
\[\sqrt{q^d}e^{itA}=H_q\]
is equivalent to
\[\sqrt{q^d}e^{i\theta_st}=e^{i\beta}\left(1+(q-1)\zeta_{6-q}\right)^{d-s}\left(1-\zeta_{6-q}\right)^s,\]
for $s=0,1,\cdots, d$. It follows that
\begin{equation}
\sqrt{q^d}e^{i\theta_0t}=e^{i\beta}(1+(q-1)\zeta_{6-q})^d\label{theta_0}
\end{equation}
and
\begin{equation}
e^{i(\theta_s-\theta_0)t}=\left(\frac{1-\zeta_{6-q}}{1+(q-1)\zeta_{6-q}}\right)^s\label{theta_s},
\end{equation}
for $s=0,1,\cdots, d$.
\begin{enumerate}[(i)]
\item For $q=2$, Equation (\ref{theta_s}) reduces to
\[\frac{2^k}{\pi}(\theta_s-\theta_0)t\equiv \epsilon 2^{k-1}s\pmod{2^{k+1}}.\]
\item For $q=3$, Equation (\ref{theta_s}) reduces to
\[\frac{3^k}{2\pi}(\theta_s-\theta_0)t\equiv\epsilon 3^{k-1}s \pmod{3^k}.\]
\item For $q=4$, Equation (\ref{theta_s}) reduces to
\[\frac{2^k}{\pi}(\theta_s-\theta_0)t\equiv 2^ks \pmod{2^{k+1}}.\]
\end{enumerate}
Thus, for $q\in\{2,3,4\}$, if $\theta_s-\theta_0$ satisfies the corresponding condition in the lemma, then there exist $t,\beta\in\re$ that satisfy Equation \ref{theta_0} and \ref{theta_s}. That is, $X$ admits uniform mixing at time $t$.
\qed

Equation \ref{ev_rec} tells us that we can check the above conditions by looking at the eigenvalues of $\cB(q-1)$.

\begin{corollary}
	Let $X$ be a graph in $\cB(q)$ with adjacency matrix
	\[a_0 + a_1A_1 + \cdots + a_d A_d.\]
	If $q=2$ and there is $\epsilon \in \{1,-1\}$ such that
	\[P^{(d-1)} \pmat{a_1\\\cdots\\a_d} \equiv \epsilon 2^{k-2} \one \pmod{2^k},\] 
	then $X$ admits uniform mixing at time $\pi/2^k$.
	If $q=3$ or $q=4$, and there is $\epsilon \in \{1,-1\}$ such that
	\[P^{(d-1)} \pmat{a_1\\\cdots\\a_d} \equiv \epsilon  q^{k-2} \one \pmod{q^{k-1}},\] 
	then $X$ admits uniform mixing at time $2\pi/q^k$.
	\qed
	\label{B(q-1)}
\end{corollary}

\section{Faster Uniform Mixing on Distance Graphs \label{distance}}
We apply the sufficient conditions developed in the last section to the distance graphs of Hamming graphs, as their eigenvalues are known. In Lemma 3.3 of \cite{Chan2013}, Chan derived a more accessible condition for uniform mixing in $\cH(d,2)$ using Property (iv) in Lemma \ref{properties}. However, this property holds only for $q=2$. To extend her result, we need more general properties for $q\in\{2,3,4\}$. The first one is a corollary to Lemma \ref{theta_diff}.

\begin{corollary}
Suppose $d\ge 1$, $r\ge 1$ and $k\ge 2$. Let $X_r$ be the $r$-distance graph of the Hamming graph $H(d,q)$, and let $\epsilon\in\{1,-1\}$. 
\begin{enumerate}[(i)]
\item If $q=2$, and
\[p_{r-1}^{(d-1)}(s)\equiv \epsilon 2^{k-2}\pmod{2^{k}},\]
for $s=0,1,\cdots,d-1$, then $X_r$ admits uniform mixing at time $\pi/2^k$.
\item If $q=3$, and 
\[p_{r-1}^{(d-1)}(s)\equiv \epsilon 3^{k-2}\pmod{3^{k-1}},\]
for $s=0,1,\cdots,d-1$, then $X_r$ admits uniform mixing at time $2\pi/3^k$.
\item If $q=4$, and
\[p_{r-1}^{(d-1)}(s)\equiv 2^{k-2} \pmod{2^{k-1}},\]
for $s=0,1,\cdots, d-1$, then $X_r$ admits uniform mixing at time $\pi/2^k$.
\end{enumerate}
\label{p_r}
\end{corollary}
\proof
We prove this for $q=2$; the other two cases are similar.

Suppose 
\[p_{r-1}^{(d-1)}(s)\equiv \epsilon 2^{k-2}\pmod{2^{k}},\]
for $s=0,1,\cdots,d-1$. By Property (iii) in Lemma \ref{properties}, this implies 
\[p_r(s+1)-p_r(s)\equiv -\epsilon 2^{k-1}\pmod{2^{k+1}}.\]
It follows that
\begin{align*}
p_r(s)-p_r(0)&=p_r(s)-p_r(s-1)+\cdots+p_r(1)-p_r(0)\\
&=-\epsilon s 2^{k-1}\pmod{2^{k+1}}.
\end{align*}
By Lemma \ref{theta_diff}, $X_r$ admits uniform mixing at time $\pi/2^k$.
\qed

With the help from the smaller scheme $\cH(d-1,q)$, we are able to construct examples in $\cH(d,q)$ that admit faster uniform mixing. It turns out that the conditions on the eigenvalues can be further simplified using Lemma \ref{properties}. From now on, we focus on the Hamming schemes $\cH(d,3)$ and $\cH(d,4)$, as the examples in $\cH(d,2)$ are already given in \cite{Chan2013}.

\begin{lemma}
For $d\ge 1$, $r\ge 1$ and $k\ge 2$, if there exists $\epsilon\in\{-1,1\}$ such that the following holds
\begin{enumerate}[(i)]
\item $2^{r-1}\binom{d-1}{r-1}\equiv \epsilon 3^{k-2}\pmod{3^{k-1}}$,
\item $3^{k-h-1}$ divides $\binom{d-h-1}{r-h-1}$ for $h=1,2,\cdots,k-2$,
\end{enumerate}
then the distance graphs $X_{r}$ and $X_{d-r+1}$ in the Hamming scheme $\cH(d,3)$ admit uniform mixing at time $2\pi/3^{k}$. \label{main3}
\end{lemma}
\proof
From Lemma \ref{properties}, we have
\begin{align*}
p_{r-1}^{(d-1)}(s)&\equiv\sum_{h=0}^{d-1}(-3)^h2^{r-h-1}\binom{d-h-1}{r-h-1}\binom{s}{h}\pmod{3^{k-1}}\\
&\equiv \sum_{h=0}^{k-2}(-3)^h2^{r-h-1}\binom{d-h-1}{r-h-1}\binom{s}{h}\pmod{3^{k-1}}.
\end{align*}
By condition (i), when $s=0$, \[p_{r-1}^{(d-1)}(0)=2^{r-1}\binom{d-1}{r-1}\equiv \epsilon 3^{k-2}\pmod{3^{k-1}}\]
and  condition (ii), when $s\ge 1$,
\begin{align*}
p_{r-1}^{(d-1)}(s)&=p_{r-1}^{(d-1)}(0)+\sum_{h=1}^{k-2}(-3)^h2^{r-h-1}\binom{d-h-1}{r-h-1}\binom{s}{h}\pmod{3^{k-1}}\\
&\equiv \epsilon 3^{k-2}\pmod{3^{k-1}}.
\end{align*}
It follows from Corollary \ref{p_r} that $X_r$ in $\cH(d,3)$ admit uniform mixing at time $2\pi/3^{k}$. For $X_{d-r+1}$, first note that condition (i) is symmetric on $r$ and $d-r+1$. By condition (ii), $3^{k-h-1}$ divides
\[\binom{d-h}{r-h}-\binom{d-h-1}{r-h-1}=\binom{d-h-1}{r-h},\]
for $h=1,2,\cdots, k-2$. It follows that $3^{k-h-2}$ divides
\[\binom{d-h-1}{r-h}-\binom{d-h-2}{r-h-1}=\binom{d-h-2}{r-h},\]
for $h=1,2,\cdots, k-2$. Continuing this procedure, we see that $3^{k-h-\ell}$ divides 
\[\binom{d-h-\ell}{r-h}\]
for $h=1,2,\cdots, k-2$ and $\ell=1,2,\cdots,k-2$. Taking $h=1$ for all $\ell$ shows that $3^{k-\ell-1}$ divides
\[\binom{d-\ell-1}{r-1}=\binom{d-\ell-1}{(d-r+1)-\ell-1},\]
for $\ell=1,2,\cdots,k-2$, which is exactly condition (ii) with $r$ replaced by $d-r+1$. Hence, $X_{d-r+1}$ admits uniform mixing at time $2\pi/3^k$ as well.
\qed

With a similar argument, we can reduce the conditions for faster uniform mixing in $\cH(d,4)$ to the following.

\begin{lemma}
For $d\ge 1$, $r\ge 1$ and $k\ge 2$, if the following two conditions hold
\begin{enumerate}[(i)]
\item $3^{r-1}\binom{d-1}{r-1}\equiv 2^{k-2}\pmod{2^{k-1}}$,
\item $2^{k-2h-1}$ divides $\binom{d-h-1}{r-h-1}$ for $h=1,2,\cdots,\lfloor k/2\rfloor-1$,
\end{enumerate}
then the distance graphs $X_r$ and $X_{d-r+1}$ in the Hamming scheme $\cH(d,4)$ admit uniform mixing at time $\pi/2^k$.
\label{main4}
\qed
\end{lemma}

The above observations imply that our potential examples rely heavily on the divisibility of a binomial coefficient by some prime power. In fact, this is closely related to the base $p$ representation of the binomial coefficients, where $p$ is prime. To find the pairs $(d,r)$ that satisfy the divisibility conditions, we need the following number theory result due to Kummer \cite[Ch 9]{Dickson1966}.

\begin{theorem}[Kummer]
Let $p$ be a prime number. The largest integer $k$ such that $p^k$ divides $\binom{N}{M}$ is the number of carries in the addition of $N-M$ and $M$ in base $p$ representation.
\qed
\end{theorem}

For our purposes, we look at the ternary representations of $d-r$ and $r-h-1$ for $\cH(d,3)$, and their binary representations for $\cH(d,4)$. The following are our new examples of distance graphs that admit uniform mixing at times earlier than the Hamming graphs.

\begin{theorem} \label{q=3(1)}
For $k\ge 2$ and $r\in\{3^k-1,3^k-4,3^k-7\}$, the $r$-distance graphs $X_r$ of the Hamming graph $H(2\cdot 3^k-9,3)$ admit uniform mixing at time $2\pi/3^{k}$.
\end{theorem}
\proof
Let $d=2\cdot 3^k-9$ and $r=3^k-1$. Then
\begin{align*}
d-r&=2\cdot 3^{k-1}+2\cdot 3^{k-2}+\cdots+2\cdot 3^2+0\cdot 3^1+1\cdot 3^0,\\
r-1-h&=2\cdot 3^{k-1}+2\cdot 3^{k-2}+\cdots+2\cdot 3^2+2\cdot 3^1+1\cdot 3^0-h.
\end{align*}
When $h=0$, since $(d-r)+(r-1)$ has exactly $k-2$ carries, $\binom{d-1}{r-1}$ is divisible by $3^{k-2}$ but not divisible by $3^{k-1}$. Then there exists $\epsilon\in\{-1,1\}$ such that 
\[2^{r-1}\binom{d-1}{r-1}\equiv \epsilon 3^{k-2}\pmod{3^{k-1}}.\]
For $h=1$, the number of carries in $(d-r)+(r-2)$ is still $2^{k-2}$. When $h=2,\cdots,k-2$, the number of carries in $(d-r)+(r-h+1)$ drops by at most one as $h$ increases by one. Therefore $(d-r)+(r-h+1)$ has at least $k-h-1$ carries, and so $3^{k-h-1}$ divides $\binom{n-h}{r-h-1}$. By Theorem \ref{main3}, $X_{3^k-1}$ and $X_{3^k-7}$ in $\cH(2\cdot 3^k-9,3)$ admit uniform mixing at time $2\pi/3^k$. Similar argument applies to $X_{3^k-4}$.
\qed

Some new examples in $\cH(d,4)$ can be obtained in a similar way.

\begin{theorem} \label{q=4}
For $k\ge 2$, the distance graph $X_{2^{k-2}}$ of the Hamming graph $H(2^{k-1}-1,4)$, and the distance graphs $X_{2^{k-1}-1}$, $X_{2^{k-1}}$ of the Hamming graph $H(2^k-2,4)$ admit uniform mixing at time $\pi/2^k$.
\qed
\end{theorem}

\section{Faster Mixing in Schemes}
In this section, we give another family of graphs in $\cH(d,q)$ that have faster uniform mixing. These are unions of some distance graphs of the Hamming graph $H(d,q)$. Compared to the examples obtained in the last section, these graphs have smaller sizes.

\begin{theorem}\label{q=3(2)}
	In the Hamming scheme $\cH(2k+1, 3)$, the graph with adjacency matrix 
	\[\sum_{\ell} A_{3\ell + i}\]
	has uniform mixing at time $2\pi/3^k$.
\end{theorem}

\proof
By Corollary \ref{B(q-1)}, it suffices to show that for each $i=0,1,2$ and the vector
\[a = \sum_{\ell} e_{3\ell +i},\]
there is $\epsilon \in \{1,-1\}$ such that
\[P^{(d-1)} a \equiv \epsilon 3^{k-1} \pmod{3^k}.\]
Also, Equation \ref{hamev} indicates that the eigenvalues $p_r^{(2k)}(s)$ for $\cH(2k,q)$ are the coefficients in 
\[f_s(x) = (1+2x)^{2k-s} (1-x)^s.\]
Let $\zeta = e^{2\pi i/3}$. We have
\begin{gather*}
f_s(1) = \begin{cases} 3^d, \quad s=0 \\ 0, \quad s\ne 0\end{cases} \\
f_s(\zeta) = -3^k \zeta^{-s}\\
f_s(\zeta^2) = -3^k \zeta^{s}.
\end{gather*}
Then
\begin{align*}
([x^0] + [x^3] + [x^6] + \cdots )f(x) 
&= \frac{1}{3} (f(1) + f(\zeta) + f(\zeta^2))\\
&= \begin{cases}
3^{2k-1} - 2\cdot 3^{k-1}, \quad s=0\\
3^{k-1} \quad s\ne 0
\end{cases}\\
&\equiv 3^{k-1} \one \pmod{3^k}.
\end{align*}
Similar computations can be carried out for 
\[ ([x^1] + [x^4] + [x^7] + \cdots )f(x) = \frac{1}{3} (f(1) + \zeta^2f(\zeta) + \zeta f(\zeta^2)),\]
and 
\[([x^2] + [x^5] + [x^8] + \cdots )f(x) = \frac{1}{3} (f(1) + \zeta f(\zeta) + \zeta^2 f(\zeta^2)). \tag*{\sqr53}\]

\section{Mixing Times}
The eigenvalues and eigenvectors of an abelian Cayley graph are determined by the group characters. For linear Cayley graphs over $\ints_q^d$, there is a simple expression of its eigenvalues in terms of the connection set $C$. With these observations, we derive a necessary and sufficient condition for uniform mixing to occur on linear Cayley graphs over $\ints_q^d$. This extends the result in \cite[Ch 5]{Zhan2014} on the case where $q=3$.
Throughout this section, the inner product $\grp{\cdot,\cdot}$ is taken over $\ints_q$.

\begin{lemma}
Let $X$ be a Cayley graph over $\ints_q^d$ with connection set $C$. For an element $a\in \ints_q^d$, let $\psi_a: \ints_q^d \to \cx$ be the map given by
\[\psi_a(x) = e^{2\pi i \grp{a,x} / q}.\]
Then $\psi_a$ is an eigenvector for $A(X)$ with eigenvalue $\psi_a(C)$. Moreover, the eigenvectors defined above are pairwise orthogonal, and they form a group isomorphic to the additive group $\ints_q^d$.  Finally, if $X$ is linear,  the eigenvalues are integers and can be computed as follows
\begin{equation}
\psi_a(C) = \frac{1}{q-1} (q \abs{C\cap a^{\perp}} - \abs{C}).
\tag*{\sqr53}
\end{equation}
\label{gpchar}
\end{lemma} 

We apply the spectral decomposition to Cayley graphs over $\ints_q^d$. Since Cayley graphs are vertex transitive, it suffices to look at the first row of the transition matrix.

\begin{lemma}
Let $g\in\ints_q^d$. The $0g$-entry of the transition matrix of $X(\ints_q^d,C)$ is
\[U_X(t)_{0,g}=\frac{1}{q^d}\sum_{a\in\ints_q^d}e^{i\psi_a(C)t}\psi_a(g).\]  \label{spec}
\end{lemma}
\proof
By Lemma \ref{gpchar}, 
\[V_\theta=\left\{\frac{1}{\sqrt{q^d}}\psi_a:\psi_a(C)=\theta\right\}\]
is an orthonormal basis of the eigenspace of $\theta$. Hence the idempotents representing the projection onto the eigenspace of $\theta$ is
\[E_{\theta}=\frac{1}{q^d}\sum_{a:\psi_a(C)=\theta}\psi_a\psi_a^*.\]
By the spectral decomposition of $U_X(t)$, we have
\[U_X(t)=\frac{1}{q^d}\sum_{a\in\ints_q^d}e^{i\psi_a(C)t}\psi_a\psi_a^*.\]
Lastly note that 
\[\psi_a(0)\comp{\psi}_a(g)=\comp{\psi}_a(g)=\psi_{-a}(g). \tag*{\sqr53}\]

In the rest of this section, we will denote the eigenvalues of $X(\ints_q^d,C)$ by
\begin{equation*}
\theta_a:=\psi_a(C) \label{theta_a}.
\end{equation*}
For linear Cayley graphs, we can characterize uniform mixing as follows.

\begin{lemma}\label{times}
Let $X$ be a linear Cayley graph over $\ints_q^d$ with connection set $C$. Uniform mixing occurs on $X$ at time $t$ if and only if for all $g\in\ints_q^d$,
\[\sum_{a,b: \grp{a-b,g}=0} e^{i(\theta_a-\theta_b)t}=q^d.\]
\end{lemma}
\proof
The condition
\[|U_X(t)_{0,g}|^2=\frac{1}{q^d} \]
is equivalent to
\begin{align}
q^d&=\left|\sum_{a\in\ints_q^d} e^{i\theta_at}e^{i 2\pi \grp{a,g}/q}\right|^2\notag \\ 
&=\sum_{a,b\in\ints_q^d}e^{i(\theta_a-\theta_b)t}e^{i2\pi\grp{a-b,g}/9}.
\label{temp}
\end{align}
Now partition pairs $(a,b)$ of group elements into $q$ classes
\[K_{\lambda} =\{(a,b): \grp{a-b, g} = \lambda\},\]
where $\lambda \in \ints_q$. Note that for any $\lambda \ne 0$, we have $(a,b) \in K_1$ if and only if $(\lambda a, \lambda b) \in K_{\lambda}$. Further, by the formula in Lemma \ref{gpchar}, 
\[\theta_a - \theta_b = \theta_{\lambda a} - \theta_{\lambda b}.\]
Therefore Equation \ref{temp} reduces to
\begin{align}
q^d &= \sum_{(a,b)\in K_0}  e^{i(\theta_a - \theta_b)t} + e^{i(\theta_a - \theta_b)t} \sum_{\lambda \ne 0} e^{i 2 \pi \lambda /q } \notag\\
&=\sum_{(a,b) \in K_0} e^{i(\theta_a - \theta_b)t} - \sum_{(a,b) \in K_1} e^{i(\theta_a - \theta_b)t}.
\label{cond1}
\end{align} 
Applying \ref{cond1} to the $00$-entry of the transition matrix, we have
\begin{equation}
q^d=\sum_{a,b}e^{i(\theta_a-\theta_b)t}=\sum_{(a,b) \in K_0}e^{i(\theta_a-\theta_b)t}+(q-1)\sum_{(a,b) \in K_1}e^{i(\theta_a-\theta_b)t}.\label{cond2}
\end{equation}
Combining \ref{cond1} and \ref{cond2} yields the desired condition.
\qed

By Lemma \ref{gpchar}, the difference between two eigenvalues of $X(\ints_q^d, C)$ is 
\[\theta_a-\theta_b=\frac{q}{q-1}\left(|C\cap a^{\perp}|-|C\cap b^{\perp}|\right)\]
which is divisible by q. Let
\[m_{ab}:=\frac{\theta_a-\theta_b}{q}.\]
We define a rational function in $x$ over the integers by
\begin{equation}
F_g(x):=\left(\sum_{a,b:\grp{a-b,g}}x^{m_{ab}}\right)-q^d. \label{f_g}
\end{equation}
Note that by symmetry in $a$ and $b$, this is a palindromic polynomial divided by some power of $x$. The mixing times of $X(\ints_q^d, C)$ are determined by the roots of these rational functions.

\begin{theorem}
$X(\ints_q^d,C)$ admits uniform mixing at time $t$ if and only if $e^{qit}$ is a zero of \[\gcd\{F_g:g\in\ints_q^d\}.\]
\label{gcd}
\end{theorem}
\proof
By Lemma \ref{times}, uniform mixing occurs at time $t$ if and only if $t$ satisfies
\[\sum_{a,b: \grp{a-b,g}=0}e^{i(\theta_a-\theta_b)t}=q^d\]
for all $g$, or equivalently, if and only if 
\[\sum_{a,b: \grp{a-b, g}=0} (e^{qit})^{\frac{\theta_a-\theta_b}{q}} - q^d=0,\]
which is exactly $F_g(e^{qit})=0$, for all $g$.
\qed

For an application of the above result, consider the linear Cayley graph $X(\ints_q^d, C)$, where $C$ consists of all non-zero multiples of $\{\seq{e}{1}{2}{d},\one\}$. It is isomorphic to the quotient graph $H(d+1,q)/\grp{\one}$. 

\begin{theorem}
If $H(d+1,q)/\grp{\one}$ admits uniform mixing at time $t$, then either 
\begin{enumerate}[(i)]
\item $q=2, 4$ and $t=k\pi/4$ for some odd $k$, or
\item $q=3$ and $t=2k \pi/9$ for some $k$ not divisible by $3$.
\end{enumerate}
\end{theorem}
\proof
We compute the function $F_{\one}(x)$. By Theorem \ref{gpchar}, the eigenvalue $\theta_a$ is determined by $\abs{C\cap a^{\perp}}$. Since the elements in $C$ are non-zero multiples of $\{\seq{e}{1}{2}{d},\one\}$, we have
\[\theta_a = 
\begin{cases}
(q-1)d - q\wt(a) + (q-1), \quad \text{ if }\grp{a, \one}=0\\
(q-1)d - q\wt(a) -1, \quad \text{ if }\grp{a,\one}\ne 0.\end{cases}\]
Now let $\alpha_j$ be the number of elements in $\grp{\one}^{\perp}$ with weight $j$. Since 
\[W_{\one}(x,y) = x^d + (q-1)y^d\]
by MacWilliams' identity,
\[W_{\one^{\perp}} = \frac{1}{q} \left( (x+(q-1)y)^d + (q-1)(x-y)^d\right).\]
Therefore 
\[n_j = \frac{1}{q} \binom{d}{j} \left((q-1)^j + (-1)^j(q-1)\right).\]
To compute the weights of the other elements in $\ints_q^d$, note that for each $\lambda \ne 0$, there is a one-to-one correspondence between $\{g: \grp{g, \one} = \lambda \}$ and $\{g: \grp{g, \one} = 1 \}$, so it suffices to compute the number of elements in 
\[\{g: \grp{g, \one} = 1 \}\]
with weight $j$, denoted $\beta_j$. Since the total number of elements in $\ints_q^d$ with weight $j$ is 
\[\binom{d}{j} (q-1)^j,\]
we have
\[(q-1) \beta_j = \binom{d}{j} (q-1)^j -\beta_j,\]
that is,
\[\beta_j = \frac{1}{q} \binom{d}{j} \left( (q-1)^j - (-1)^j)\right).\] 
Hence,
\begin{align*}
F_{\one}(x) &=\sum_{a,b} x^{m_{ab}} - q^d\\
&=\sum_{a, b: \grp{a,\one} = \grp{b,\one}=1} x^{\wt(a)-\wt(b)} + (q-1) \sum_{a,b:\grp{a,\one} = \grp{b,\one}=1} x^{\wt(a)-\wt(b)} - q^d\\
&=\sum_j\sum_k (\alpha_j \alpha_k+(q-1) \beta_j\beta_k) x^{j-k} -q^d\\
&=\frac{1}{q}\sum_j \sum_k\binom{d}{j}\binom{d}{k}\left((q-1)^{j+k} + (-1)^{j+k}(q-1)\right) x^{j-k} -q^d\\
&=\frac{1}{q}\left((q-1)\left(x+\frac{1}{x}\right) + (q-1)^2+1\right) + \frac{q-1}{q} \left(2-\left(x+\frac{1}{x}\right)\right)^d-q^d.
\end{align*}
Now let $z=x+1/x$. Recall that $x=e^{iqt}$ for some $t$, so $-2\le z \le 2$. Substitute $z$ into $F_{\one} (x)$ and we have
\[F_{\one}(z):= \frac{1}{q} ((q-1)z + (q-1)^2 +1)^d+ \frac{q-1}{q} (2-z)^d - q^d.\]
The derivative of $F_{\one}(z)$ is positive if and only if 
\[((q-1)z + (q-1)^2+1)^{d-1} >  (2-z)^{d-1}.\]
Notice that the expression in the brackets of the left hand side is at least
\[-2(q-1)+(q-1)^2+1=(q-2)^2\ge 0,\]
so $F_{\one}'(z)>0$ if and only if 
\[(q-1)z + (q-1)^2+1 > 2-z,\]
that is, $z> 2-q$. It follows that for $q>4$, 
\[F_{\one}(z) \le F_{\one}(2) = q^{d-1}-q^d<0.\]
Hence $H(d,q)/\grp{\one}$ does not admit uniform mixing when $q>4$. For $q\in\{2,3,4\}$, we see that $z=2-q$ is the stationary point and the only zero in the interval $[-2,2]$. Therefore uniform mixing must occur at time $t$ for which $2\cos(qit) = 2-q$.
\qed

For all the Cayley graphs over $\ints_q$ known to admit uniform mixing, the mixing time is of the form $2\pi/qn$ for some integer $n$. As a second consequence of Theorem \ref{gcd}, the degree of such a graph must be large enough for mixing to occur at time $2\pi/qn$.

\begin{corollary}
Let $\phi(n)$ be the Euler's totient function. If uniform mixing occurs on $X(\ints_q^d, C)$ at time $2\pi/qn$, then 
\[\abs{C} \ge \frac{q-1}{2} (\phi(n)+q-1).\]
\label{totient}
\end{corollary}

\proof
Let $g\in C$. For $a,b\in \ints_q^d$ such that $\grp{a-b,g}=0$,
\[|C\cap a^{\perp}|-|C\cap b^{\perp}|\le |C|-2.\]
Hence 
\[\deg(f_g)\le |C|-q+1.\]
If $X(\ints_q^d,C)$ admits uniform mixing at $2\pi/qn$, then $F_g$ is divisible by the $n$-th cyclotomic polynomial $\Phi_n(n)$ with degree $\phi(n)$. Thus 
\[\phi(n)\le \frac{2}{q-1}(|C|-q+1).\tag*{\sqr53}\]

\section{Local and Global Uniform Mixing on Stars \label{star}}

For irregular graphs, local uniform mixing may be a better choice to start with.  We follow Carlson et al \cite{Carlson2006} and show that the star $K_{1,n}$ admits local uniform mixing. In particular, uniform mixing in the global sense occurs on the claw $K_{1,3}$.

We apply spectral decomposition to the adjacency matrix $A$ of the star $K_{1,n}$. The eigenvalues of $K_{1,n}$ are $\theta_0=0$, $\theta_1=\sqrt{n}$ and $\theta_2=-\sqrt{n}$. Denote the projections onto these eigenspaces by $E_0$, $E_1$ and $E_2$. We have
\begin{gather*}
E_0=\pmat{0 & \zero^T\\\zero & I-\frac{1}{n}J},\\
E_1=\frac{1}{2n}\pmat{n & \sqrt{n}\one^T\\\sqrt{n}\one & J},\\
E_2=\frac{1}{2n}\pmat{n & -\sqrt{n}\one^T\\-\sqrt{n}\one & J},
\end{gather*}
where $J$ denotes the all-ones matrix. It follows that the transition matrix of $K_{1,n}$ is 
\begin{align*}
U(t)&=e^{0\cdot it}E_0+e^{\sqrt{n}it}E_1+e^{-\sqrt{n}it}E_2\\
&=\pmat{\cos\left(\sqrt{n}t\right) & \frac{i}{\sqrt{n}}\sin\left(\sqrt{n}t\right)\one\\ \frac{i}{\sqrt{n}}\sin\left(\sqrt{n}t\right)\one & I+\frac{1}{n}\left(\cos\left(\sqrt{n}t\right)-1\right)J}.
\end{align*}
The quantum walk starting with the central vertex is uniform mixing at time $t$ if and only if 
\[\left|\cos\left(\sqrt{n}t\right)\right|=\left|\frac{\sin\left(\sqrt{n}t\right)}{\sqrt{n}}\right|\]
or equivalently, 
\begin{equation}\tan\left(\sqrt{n}t\right)=\pm \sqrt{n}.
\label{local}
\end{equation}
Thus, the star $K_{1,n}$ admits local uniform mixing at time 
\[\pm \frac{\arctan\left({\sqrt{n}}\right)}{\sqrt{n}}+k\pi\] for all integers $k$.

For uniform mixing, one additional condition from the lower right block of $U(t)$ is 
\[\left|1+\frac{1}{n}\left(\cos\left(\sqrt{n}t\right)-1\right)\right|=\left|\frac{1}{n}\left(\cos\left(\sqrt{n}t\right)-1\right)\right|\]
or equivalently,
\begin{equation}
\cos\left(\sqrt{n}t\right)=1-\frac{n}{2}.
\label{global}
\end{equation}
Combining Equation (\ref{local}) and Equation (\ref{global}), we see that the only solution is 
\[n=3,\quad t=\pm \frac{2\pi}{\sqrt{27}}+2k\pi\]
for all integers $k$. Plugging this into $U(t)$ yields a flat matrix. We conclude that the only star that admits uniform mixing is the claw $K_{1,3}$, with earliest mixing time $2\pi/\sqrt{27}$. The Cartesian powers of $K_{1,3}$ then form an infinite family of irregular graphs that admit uniform mixing.

\section{Open Problems}

There are a number of open problems on unifom mixing, ranging across characterizing graphs that admit uniform mixing in some common family, determining the mixing times of a given graph, and constructing new examples. Following our notation in Section \ref{quotient}, we let $\tau_q$ denote the earliest time at which the complete graph $K_q$ admits uniform mixing.

\begin{enumerate}
\item Question: To determine whether uniform mixing occurs on the quotient graph $H(d,q)/\Gamma$ at time $\tau_q$, we have to check the weight distribution of every coset of $\Gamma$. Is it sufficient to just check the weight distribution of $\Gamma$?

For groups with one or two generators, the weight distribution of any coset $\Gamma+c$ is merely a permutation of the weight distribution of $\Gamma$. For an example see Theorem \ref{2-gen}. If this were true in general, it would be helpful in characterizing linear Cayley graphs with higher degrees.

\item Question: Is there a characterization of uniform mixing on non-linear Cayley graphs over $\ints_q^d$?

This is one thing that the weight distribution condition does not tell. For $q\ge 4$, a Cayley graph over $\ints_q^d$ may not be linear, and thus may not be a quotient graph of $H(d,q)$. It is desirable to find another approach for these non-linear Cayley graphs.

\item Question: If a Cayley graph over $\ints_q^d$ admits uniform mixing, must its eigenvalues be integral?

As we mentioned in Section \ref{intro}, all the known Cayley graphs that admit uniform mixing have integer eigenvalues. It is unclear if this is a necessary condition. The first place to find a counterexample might be the Cayley graphs over $\ints_5^d$, as the eigenvalues are no longer guaranteed to be integral. 

\item Question: If a Cayley graph over $\ints_q^d$ admits uniform mixing at time $t$, must $t$ be a rational multiple of $\pi$?

Again this is true for all the known examples. However, it may only apply to graphs with integer eigenvalues. Even for this smaller class of graphs, it would  be interesting to confirm such an algebraic property of the mixing times.

\item Question: How fast can a Cayley graph over $\ints_q^d$ admit uniform mixing?

So far, the best examples that admit  uniform mixing earlier than $\tau_q$ are the distance graphs of $H(d,2)$ found in \cite{Chan2013}, and the distance graphs of $H(d,3)$ and $H(d,4)$ found in Section \ref{distance} of this paper. These families provide arbitrarily faster uniform mixing, although at the cost of larger vertex sets. In an effort to construct new examples with faster uniform mixing, a question arises as to whether there is a lower bound on the mixing time of a given graph.

\item Question: Are there more irregular graphs that admit uniform mixing?

The star $K_{1,3}$ and its Cartesian powers suggest that there could be other irregular graphs that admit uniform mixing. As we did in Section \ref{star}, one may look at local uniform mixing on some common families of irregular graphs, and then impose more conditions for global uniform mixing.

\end{enumerate}

\bibliographystyle{amsplain}
\bibliography{um.bib}

\end{document}